\documentclass[12pt]{article}

\RequirePackage{kpfonts}
\RequirePackage[sf,bf,small,raggedright]{titlesec}
\usepackage{textcomp}
\usepackage{amssymb}
\usepackage{bbm}
\usepackage{fancyhdr}
\usepackage{amsmath}
\usepackage{amsbsy,amsthm}
\usepackage{amscd}
\usepackage{latexsym}
\usepackage{graphicx}   
\usepackage{pdfsync}
\usepackage{blkarray}
\usepackage{multirow}
\usepackage{hyperref}
\usepackage{color}
\usepackage{comment}
\usepackage{framed}

\newcommand{\R}{\mathbb{R}}

\newcommand{\A}{\mathcal{A}}

\newtheorem{lemma}{Lemma}
\newtheorem{theorem}{Theorem}

\newtheorem{proposition}{Proposition}

\numberwithin{equation}{section}

\newcommand{\argmin}{\mathop{\mathrm{argmin}}}

\newcommand{\EXP}{\mathbb{E}}

\newcommand{\PROB}{\mathbb{P}}

\begin{document}

\title{Finding the seed of uniform attachment trees
\thanks{
G\'abor Lugosi was supported by
the Spanish Ministry of Economy and Competitiveness,
Grant MTM2015-67304-P and FEDER, EU.}
}
\author{
G\'abor Lugosi\thanks{Department of Economics and Business, Pompeu
  Fabra University, Barcelona, Spain, gabor.lugosi@upf.edu}
\thanks{ICREA, Pg. Lluís Companys 23, 08010 Barcelona, Spain}
\thanks{Barcelona Graduate School of Economics}
\and
Alan S. Pereira\thanks{Instituto Nacional de Matemática Pura e Aplicada (IMPA), Estrada Dona Castorina, 110 - Jardim Bot\^anico, Rio de Janeiro - RJ, Brazil, alanand@impa.br}
}

\maketitle

\begin{abstract}
A uniform attachment tree is a random tree that is generated dynamically. Starting from a fixed ``seed'' tree,
vertices are added sequentially by attaching each vertex to an existing vertex chosen uniformly at random.
Upon observing a large (unlabeled) tree, one wishes to find the initial seed.
We investigate to what extent seed trees can be recovered, at least partially. We consider 
three types of seeds: a path, a star, and a random uniform attachment tree. We propose and analyze
seed-finding algorithms for all three types of seed trees.
\end{abstract}

\section{Introduction}

Dynamically growing networks represent complex relationships in numerous areas of science. 
In a rapidly increasing number of applications, one does not observe the entire dynamical growth
procedure but merely a present-day snapshot of the network is available for observation. Based on this snapshot,
one wishes to infer various properties of the \emph{past} of the network. Such problems belong to 
the area that may be termed \emph{network archeology}, see Navlakha and Kingsford \cite{NaKi11}.

The simplest dynamically grown networks are trees that are grown by attaching vertices sequentially
to the existing tree at random, according to a certain rule. In the \emph{uniform attachment} model,
at each step, an existing vertex is selected uniformly ar random, and a new vertex is attached to it 
by an edge. When the process is initialized from a single vertex, this procedure gives rise to the well-studied
\emph{uniform random recursive tree}, see Drmota \cite{Drm09}. In \emph{preferential attachment} models
(such as plane-oriented recursive trees) existing vertices with higher degrees are more likely to be chosen
to be attached to. In this paper we consider randomly growing uniform attachment trees that are grown from
a fixed \emph{seed}. Thus, initially, the tree is a given fixed (small) tree and further vertices are attached
according tio the uniform attachment process.

``Archeology'' of randomly growing trees has received increasing attention recently, see
Brautbar and Kearns \cite{BrKe10},
Borgs, Brautbar, Chayes, Khanna, and Lucier \cite{BoBrChKhBr12},
Bubeck, Devroye, and Lugosi \cite{BuDeLu17},
Bubeck, Mossel, and R{\'a}cz \cite{BuMoRa15},
Bubeck, Eldan Mossel, and R{\'a}cz \cite{BuElMoRa17},
Curien, Duquesne, Kortchemski, and Manolescu \cite{CuDuKoMa15},
Frieze and Pegden \cite{FrPe14},
Jog and Loh \cite{JoLo17a,JoLo17b},
Shah and Zaman \cite{ShZa11,ShZa16}
for a sample of the growing literature.

Several papers consider the problem of finding the initial vertex (or root) in a randomly growing tree started from
a single vertex, see 
Brautbar and Kearns \cite{BrKe10},
Borgs, Brautbar, Chayes, Khanna, and Lucier \cite{BoBrChKhBr12},
Frieze and Pegden \cite{FrPe14},
Shah and Zaman \cite{ShZa11,ShZa16}, 
Bubeck, Devroye, and Lugosi \cite{BuDeLu17},
Jog and Loh \cite{JoLo17a,JoLo17b} for various models. Randomly growing trees started from an 
initial seed tree were considered by 
Bubeck, Mossel, and R{\'a}cz \cite{BuMoRa15},
Bubeck, Eldan Mossel, and R{\'a}cz \cite{BuElMoRa17}, and
Curien, Duquesne, Kortchemski, and Manolescu \cite{CuDuKoMa15}. These papers prove that in uniform
and preferential attachment models, for any pair of possible seed trees, one may construct a hypothesis test
that decides which of the two seeds generated the observed tree, with a probability of error strictly smaller than $1/2$,
regardless of the size of the observed tree.

In this paper we consider the problem of \emph{finding} the seed tree (of known structure) in a large observed tree. 
The questions we seek to answer are: (1) to what extent is it possible to identify the seed tree?
(2) what is the role of the structure of the seed in the difficulty of the reconstruction problem?
While we are far from completely answering these questions, this paper contributes to the understanding
of these problems. In particular, we consider three types of possible seed trees, namely paths, stars, and 
random uniform recursive trees. For each of these examples, we present algorithms to recover, at least 
partially, the seed tree. In all cases, partial recovery is possible, with any prescribed probability of error,
regardless of the size of the observed tree. However, the difficulty of the recovery depends heavily on the 
structure of the tree. Paths and stars are considerably easier to find than uniform random recursive trees.

In Section \ref{sec:setup} we introduce the mathematical model and state the main
results. The proofs of all results are presented in Section \ref{sec:proofs}.

\section{Setup and results} 
\label{sec:setup}

Let $\ell\ge 1$ be a positive integer and let $S_\ell$ be a tree (i.e., a connected acyclic graph) 
on the vertex set $\{1,\ldots,\ell\}$. Let $n > \ell$ be another positive integer.
We say that a random tree $T_n$ on the vertex set $\{1,\ldots,n\}$ is a 
\emph{uniform attachment tree with seed} $S_\ell$ if it is generated as follows:

\begin{enumerate}
\item 
$T_{\ell}= S_{\ell}$~;
\item
For $\ell< i \le n$, $T_i$ is obtained from $T_{i-1}$ by joining vertex $i$ to a vertex of $T_{i-1}$
chosen uniformly at random, independently of all previous choices.
\end{enumerate}

The problem we study in this paper is the following. Suppose one observes a tree $T_n$ generated
by the uniform attachment process with seed $S_\ell$ but with the vertex labels hidden. The goal is
to find the seed tree $S_\ell$ in the observed unlabeled tree. More precisely, given a target accuracy
$\epsilon \in (0,1)$ a seed-finding algorithm of \emph{first kind}
outputs a set $H_1(T_n,\epsilon)$ of vertices of size $k_{\ell}\le \ell$, such that, 
with probability at least $1-\epsilon$, 
$H_1(T_n,\epsilon) \subset S_\ell$, that is, all elements of $H_1(T_n,\epsilon)$ are vertices of the
seed tree $S_\ell$. (Here, with a slight abuse of notation, we identify the seed $S_\ell$ with its vertex set $\{1,\ldots,\ell\}$.)

Similarly, a seed-finding algorithm of \emph{second kind}
outputs a set $H_2(T_n,\epsilon)$ of vertices of size $k_{\ell}\ge \ell$, such that, 
with probability at least $1-\epsilon$, 
$S_\ell \subset H_2(T_n,\epsilon)$, that is,  $H_2(T_n,\epsilon)$ contains all vertices of the
seed tree $S_\ell$. 

In both cases, one would like to have $k_\ell$ as close to $\ell$ as possible, 
even for small values of $\epsilon$. 

Bubeck, Devroye, and Lugosi \cite{BuDeLu17} considered the case $\ell=1$, that is, when the seed tree is a 
single vertex and seed-finding algorithms of the second kind. Thus, the aim of the seed-finding algorithm 
is to find the root of the observed tree.
Their main finding is that, for all $\epsilon$, the optimal value
of $k_1$ stays bounded as the size $n$ of the observed tree goes to infinity. They also show that there
exist seed-finding algorithms of the second kind such that $k_1 = o(\epsilon^{-a})$ for all $a>0$.

In this paper we show that, if $\ell$ is sufficiently large (depending on $\epsilon$), then $k_\ell$ may be
made \emph{proportional} to $\ell$ for  seed-finding algorithms of second kind, and we make similar
statements for $k_\ell$  for certain seed-finding algorithms of first kind.
How  the required value of $\ell$ depends on $\epsilon$ and what the achievable proportions are depend
heavily on the structure of the seed. We consider three prototypical examples of seeds:

$\bullet$
A \emph{path} $P_\ell$ on $\ell$ vertices is a tree that has exactly two vertices 
of degree one and $\ell-2$ vertices of degree two.

$\bullet$
A \emph{star} $E_\ell$ on $\ell$ vertices is a tree that has $\ell-1$ vertices of degree one and one vertex of degree $\ell-1$.

$\bullet$
The third example we consider is when the seed $S_\ell$ is a uniform random recursive tree on $\ell$ vertices.
In this case the proposed seed finding algorithm does not need to know the structure of the tree. Thus, this 
example may be considered as a generalization of the root-finding problem studied in \cite{BuDeLu17}. Here, 
instead of trying to locate the root of the tree, the goal is to find the first $\ell$ generations of the observed   
uniform random recursive tree $T_n$.

In what follows we present the main findings of the paper that establish the existence of seed-finding 
algorithms that are able to recover a constant fraction of the seed if it is a uniform random recursive tree.
If the seed is either a path or a star, then the situation is even better as one can recover almost the entire seed.

Importantly, all bounds established below are independent of the size $n$ of the observed tree, meaning
that (partial) reconstruction of the seed is possible regardless of how large the observed tree $T_n$ is.

\subsection{Finding the seed when it is a path}

We begin with the case when the seed is a path:

\begin{theorem}
\label{thm:path}
Let  $\epsilon \in (0,1)$ and $\gamma \in (0,1)$ and
let  $\ell \geq \max\left\{ \dfrac{2e^2}{\gamma} \log  \dfrac{1}{\epsilon}  , \dfrac{2e^2}{\gamma}\log (4e^2)\right\}$
be a positive integer.  Then for all $n\ge \ell$ sufficiently large, if
 $T_n$ is a uniform attachment tree with seed $S_\ell=P_\ell$ (a path of $\ell$ vertices), then
there exists a seed-finding algorithm that outputs a vertex set $H_n\subset \{1,\ldots,n\}$ 
with $|H_n| \geq (1-\gamma)\ell$ such that
\[
 \PROB\left\{ H_n \subset P_\ell \right\} \geq 1-\epsilon~.
\]
\end{theorem}

The theorem states that, for any fixed $\gamma>0$, if the size of the seed path $\ell$ is at least of the order of $\log(1/\epsilon)$,
then there exists an algorithm that finds all but a $\gamma$-fraction of the seed path, regardless of how large the 
observed tree $T_n$ is. Note that the required length of the path is merely logarithmic in $1/\epsilon$.
In fact, this dependence is essentially best possible. The following result shows that if the seed path has less than
$ \frac{ \log(1/\epsilon)}{\log \log (1/\epsilon)}$ vertices, then \emph{any} seed finding algorithm must miss 
at least half of the seed, with probability greater than $\epsilon$.

\begin{theorem}
\label{thm:pathlower}
Let $\epsilon\in (0,e^{-e^2})$. 
Suppose that $T_n$ is a uniform attachment tree with seed $S_\ell=P_\ell$ for 
$\ell \leq  \frac{ \log(1/\epsilon)}{\log \log (1/\epsilon)}$. Then, for all $n\ge 2\ell$,
any seed-finding algorithm that outputs a vertex set $H_n$ of size $\ell$
has
\[
 \PROB\left\{ |H_n \cap P_\ell| \le \frac{\ell}{2} \right\} \geq \epsilon~.
\]
\end{theorem}

\subsection{Finding the seed when it is a star}

Next we state our results for the case when the seed tree is a star $E_\ell$ on $\ell$ vertices. 

\begin{theorem}
\label{thm:star}
There exists a numerical positive constant $C$ such that the following holds.
Let  $\epsilon \in (0,1)$ and $\gamma \in (0,1)$ and
let  $\ell \geq \max(C,8/\gamma)\log(1/\epsilon)$
be a positive integer.  Then for all $n\ge \ell$ sufficiently large, if
 $T_n$ is a uniform attachment tree with seed $S_\ell=E_\ell$ (a star of $\ell$ vertices), then
there exists a seed-finding algorithm that outputs a vertex set $H_n\subset \{1,\ldots,n\}$ 
with $|H_n| \leq (1+\gamma)\ell$ such that
\[
 \PROB\left\{ E_\ell \subset H_n \right\} \geq 1-\epsilon~.
\]
\end{theorem}

Once again, the order of magnitude for the required size of the seed
star is essentially optimal as a function of $\epsilon$. The proof of the
next theorem is similar to that of Theorem \ref{thm:pathlower} and
thus it is omitted.

\begin{theorem}
\label{thm:starlower}
Let $\epsilon\in (0,e^{-e^2})$. 
Suppose that $T_n$ is a uniform attachment tree with seed $S_\ell=E_\ell$ for 
$\ell \leq  \frac{ \log(1/\epsilon)}{\log \log (1/\epsilon)}$. Then, for all $n\ge 2\ell$,
any seed-finding algorithm that outputs a vertex set $H_n$ of size $\ell$
has
\[
 \PROB\left\{ |H_n \cap E_\ell| \le \frac{\ell}{2} \right\} \geq \epsilon~.
\]
\end{theorem}

\subsection{Finding the first generations}

Finally, we consider the case when the seed tree is a uniform random recursive
tree in $\ell$ vertices. Unlike in the previous two examples, here the
seed finding algorithm does now ``know'' the exact structure of the
seed. This model may be equivalently formulated as follows: starting from
a single vertex, one grows
a uniform random recursive tree $T_n$ of $n$ vertices. Upon observing
$T_n$ (without vertex labels), one's aim is to recover as much of the
tree $T_\ell$ (containing vertices attached in the first $\ell$
generations) as possible. The next theorem establishes the existence of
a seed-finding algorithm of the first kind that identifies an $\Omega(1/\log(1/\epsilon))$ fraction
of the vertices of the seed $T_\ell$ with probability at least $1-\epsilon$,
whenever $\ell$ is at least proportional to
$\log^3(1/\epsilon)$.
One should note that this result is weaker than the one obtained
for seed paths and seed stars above in various ways.
First, unlike in the cases of Theorems \ref{thm:path} and \ref{thm:star},
here we cannot guarantee that almost all of the seed tree is identified,
but only a fraction of it whose size depends on $\epsilon$--although
in a mild manner. Second, the size of the
seed tree needs to be somewhat larger as a function of $\epsilon$
as before. While in the previous cases $\ell$ needed to be logarithmic
in $1/\epsilon$, now it needs to scale as $\log^3(1/\epsilon)$.
Below we show that to some extent
these weaker results are inevitable and that finding the seed tree $T_\ell$
is inherently harder than finding more structured seed trees such as
stars and paths.

Our main positive result is as follows.

\begin{theorem}
\label{thm:oldtestament}
Let $T_n$ be a uniform random recursive tree on $n$ vertices and let
$\epsilon>0$ and $\ell\ge 1$.
Let $a=2\log(4\ell^2/\epsilon)+1$.
If $\ell$ is so large that
\[
\ell  \ge 64a^2\log(22a\ell^2/\epsilon)~,
\]
then
there exists a seed-finding algorithm that outputs a vertex set $H_n\subset \{1,\ldots,n\}$ with $|H_n| \geq \ell/(3a)$ such that
\[
\liminf_{n \to \infty} \PROB\left\{ H_n \subset T_{\ell} \right\} \geq 1-\epsilon~.
\]
\end{theorem}

Note that the condition for $\ell$ is satisfied for $\ell \ge
C\log^2(1/\epsilon)$ for a constant $C$.

Next we show that, regardless how large $\ell$ is, for $n$
sufficiently large any seed-finding algorithm of first kind
needs to output a set of vertices whose size is at most $c\ell$
where $c$ is strictly smaller than $1$. 
Similarly, any seed-finding algorithm of second kind
needs to output a set of vertices whose size is at least $C\ell$
where $C>1$. 

 In other words, 
when the seed tree is a uniform random recursive tree, the problem
of finding it is strictly harder than finding a seed path or a seed star
in the sense that no algorithm 
can have a performance as the one established in 
Theorem \ref{thm:path} or Theorem \ref{thm:star}. 
Note however, that there remains a gap between
the performance bound of Theorem \ref{thm:oldtestament} and
the impossibility bound of Theorem \ref{thm:oldtestamentlower} below,
as
the size of the vertex set in the seed found by the algorithm of
Theorem \ref{thm:oldtestament} is only guaranteed to be of the order
of $\ell/\log(1/\epsilon)$, a linear fraction but depending on
$\epsilon$. 

The impossibility results mentioned above follow from the fact that,
at time $2\ell$, a linear fraction of the vertices of the seed $T_\ell$ 
become indistinguishable from vertices that arrive between time
$\ell+1$ and $2\ell$. To make the statement precise, we need a few
definitions.

In a uniform random recursive tree $T_\ell$, we call a vertex a
\emph{singleton}
if it is a leaf and it is the only descendant of its parent vertex.

Now consider a vertex $v$ in $T_\ell$ and its position in the tree $T_{2\ell}$.
We say that $v$ is a \emph{camouflaging} vertex if

\begin{enumerate}
\item 
In $T_\ell$, $v$ is a parent of a singleton $d$;
\item
Between time $\ell+1$ and $2\ell$ a vertex $w$ is attached to $v$ such
that $w$ is a leaf of $T_{2\ell}$
\item
$d$ is a leaf of $T_{2\ell}$.
\end{enumerate}

Clearly, at time $2\ell$, and therefore at any time $n\ge 2\ell$,
the two descendants $d$ and $w$ of any camouflaging vertex $v$ are
indistinguishable. Let $G_{\ell}$ denote the number of camouflaging vertices.
Then if a seed-finding algorithm outputs 
a vertex set that contains an $(1-\gamma)\ell$ vertices of the seed,
then one must have $G_{\ell} < \gamma \ell$. The next proposition 
shows that $\gamma \ge 1/384$ with high probability.

\begin{theorem}
\label{thm:oldtestamentlower}
  For any $\ell\ge 1$,
\[
   \EXP G_\ell \ge \frac{\ell}{384}
   \]
and for any $t\geq 0$,
\[
\PROB\left\{G_\ell \leq  \frac{\ell}{384} - t \right\} \leq e^{\frac{-t^2}{2\ell}}~.
\]
\end{theorem}

\section{Proofs}
\label{sec:proofs}

In this section we present the proofs of all theorems. The construction of all seed-finding algorithms
uses a simple notion of centrality that we recall first.

\subsection{Centrality}
\label{sec:centrality}

Let $T$ be a tree with vertex set $V(T)$. 
A \emph{rooted tree} $(T,v)$ is the tree $T$ with a distinguished vertex $v\in V(T)$.
For a vertex $u\in V(T)$, denote by
 $(T,v)_{u\downarrow}$ the rooted subtree of $T$ whose root is $u$ and whose vertex set
contains all vertices $w$ of $V(T)$ such that the (unique) path connecting $w$ and $v$ in $T$  
contains $u$.

Given tree $T$, the \emph{anti-centrality} of a vertex $v \in V(T)$ is defined by 
\[
\psi(v)= \max_{u \in V(T) \setminus \{v\}} \left|(T,v)_{u\downarrow}\right|~.
\]
Thus, $\psi(v)$ is the size of the largest subtree of the tree $T$ rooted at $v$. Note that leaves of a tree $T$ have the largest anti-centrality with
$\psi(v)=|V(T)|-1$.
We say that $v$ is \emph{at least as central as} $w$ if $\psi(v)\leq \psi(w)$.

For a positive integer $k$, we denote by $H_\psi(k)$ the set of $k$ vertices of
with smallest anti-centrality, where ties may be broken arbitrarily.

This notion of centrality played a crucial role in some of the root-finding algorithms of \cite{BuDeLu17}.
We refer to Jog and Loh \cite{JoLo17a,JoLo17b} for a study of this notion in various random tree models,
including uniform random recursive trees.

\subsection{Proof of Theorem \ref{thm:path}}

Let $\epsilon,\gamma$, and $\ell$ be as in the assumptions of the theorem. 
We may assume, without loss of generality, that $\gamma\ell/2$ is an integer.
We analyze a simple seed-finding algorithm that achieves the performance stated in the theorem.
The proposed algorithm simply takes the $(1-\gamma)\ell$ most central vertices, as measured by
the function $\psi$ defined in Section \ref{sec:centrality}. 

Formally, let $k_\ell=(1-\gamma) \ell$ and define $H_n = H_{\psi}(k_\ell)$ be the set
of $k_\ell$ most central vertices of the observed tree $T_n$.

It suffices to prove that, for all sufficiently large $n$, with probability at least $1-\epsilon$, all
vertices of $T_n$ not in the seed $P_\ell$ are less central than any vertex in $P_\ell$ whose distance to 
the leaves of $P_\ell$ is at least $\gamma \ell/2$, that is,
\begin{equation}
\label{eq:pathcentr}
   \PROB\left\{ \min_{\ell<i\le n} \psi(i) > \max_{\ell \gamma/2 \le j \le \ell (1-\gamma/2)} \psi(j) \right\} \ge 1-\epsilon~.
\end{equation}
(Recall that the vertex set of the seed $P_\ell$ is $\{1,\ldots,\ell\}$.)

Let $C_1,\ldots,C_\ell$ denote the components of the forest obtained by removing the edges of $P_\ell$ from $T_n$
such that $k\in C_k$ for $k=1,\ldots,\ell$. Then
\begin{eqnarray*}
\PROB\left\{ \min_{\ell<i\le n} \psi(i)\le \max_{\ell \gamma/2 \le j \le \ell (1-\gamma/2)} \psi(j) \right\} 
& \le  &
\sum_{j=\gamma\ell/2}^{(1-\gamma/2)\ell} \PROB\left\{ \min_{\ell<i\le n} \psi(i) \leq \psi(j) \right\} \\
& \le & 
\sum_{j=\gamma\ell/2}^{(1-\gamma/2)\ell} \sum_{k=1}^{\ell} \PROB\left\{ \exists v \in C_k\setminus \{k\}: \psi(v) \leq \psi(j) \right\}~.
\end{eqnarray*}
To bound the probabilities on the right-hand side, 
suppose, without loss of generality, that $k\leq j$.  (The case $k>j$ is analogous.) 
If $v \in C_k\setminus \{k\}$ is such that $\psi(v) \leq \psi(j)$. Let $u$ be a vertex connected to $v$ 
such that $\left|(T,v)_{u\downarrow}\right|$ is maximal (i.e., $\psi(v)= \left|(T,v)_{u\downarrow}\right|$).
Then there are two possibilities:

\begin{itemize}
\item[(a)]
$(T,v)_{u\downarrow}$ is contained in $C_k$. In this case $|C_k| \geq \sum_{i \neq k}|C_i|$; 
\item[(b)]
$(T,v)_{u\downarrow} =\left(\bigcup_{i=1, i\neq k}^{\ell} C_i \right) \cup C_k'$ for some $C_k' \subset C_k$.
In this case  
\[
\left|\bigcup_{i\neq k} C_i\right| \leq \psi(v) \leq \psi(j) \leq \left|\bigcup_{i=1}^j C_i\right|
\]
which implies $\sum_{i=j+1}^{\ell} |C_i| \leq  |C_k|$.
\end{itemize}

By this observation, we have
\begin{eqnarray*}
\PROB\left\{ \exists v \in C_k\setminus \{k\}: \psi(v) \leq \psi(j) \right\}
&\leq & \PROB\left\{|C_k| \geq \sum_{i \neq k}|C_i|\right\} 
+ \PROB \left\{ \sum_{i=j+1}^{\ell} |C_i| \leq  |C_k| \right\}\\
& \leq & \PROB\left\{ |C_k| \geq \sum_{i \neq k}|C_i|\right\} 
+ \PROB \left\{ \sum_{i=(1-\gamma/2)\ell}^{\ell} |C_i| \leq  |C_k| \right\} 
\end{eqnarray*}
Now let $t=\gamma/e^2$. Then the right-hand side of the inequality above may be 
bounded further by
\begin{eqnarray*}
\PROB\left\{\sum_{i=1,i\neq k}^{\ell} |C_i| \leq nt \right\} +  \PROB\left\{ \sum_{i=1}^{\gamma \ell} |C_i| \leq nt\right\} + 2\PROB\left\{  |C_k| \ge nt \right\}
\end{eqnarray*}
Thus, we have
\begin{eqnarray*}
\lefteqn{
\PROB\left\{ \min_{\ell<i\le n} \psi(i)\le \max_{\ell \gamma/2 \le j \le \ell (1-\gamma/2)} \psi(j) \right\} 
}\\
& \le  &
(1-\gamma)\ell^2 \left( \PROB\left\{\sum_{i=1,i\neq k}^{\ell} |C_i| \leq nt \right\}  +  \PROB\left\{\sum_{i=1}^{\gamma \ell} |C_i| \leq nt\right\}+ 2\PROB\left\{ |C_k| \ge nt \right\}
 \right)
\end{eqnarray*}
To understand the behavior of the probabilities on the right-hand side, note that, for any $k=1,\ldots,\ell-1$,
$\sum_{i=1}^{k} |C_i|$ is just the number of red balls after taking $n$ samples in a standard P\'olya urn initialized with
$k$ red and $\ell-k$ blue balls. This implies that $\sum_{i=1}^{k} |C_i|/n$
converges, in distribution, to a $\text{Beta}(k,\ell-k)$ random variable.
Hence,
\[
   \lim_{n\to\infty}\PROB\left\{ |C_k|/n \ge t \right\} = (1-t)^{\ell-1} \le e^{-t(\ell-1)}
\]
and
\begin{eqnarray*}
\lim_{n\to\infty} \PROB\left\{\sum_{i=1,i\neq k}^{\ell} |C_i|/n \leq t \right\}  & \le &
\lim_{n\to\infty}  \PROB\left\{\sum_{i=1}^{\gamma \ell} |C_i|/n \leq t\right\} \\
& = & (\ell-1) \binom{\ell-1}{\gamma\ell-1} \int_0^t x^{\gamma\ell-1}(1-x)^{\ell-\gamma\ell-1} dx~.
\end{eqnarray*}
We may bound the expression on the right-hand side by
\[
   \frac{\ell^{\gamma\ell}}{(\gamma\ell-1)!}  \int_0^t x^{\gamma\ell-1} dx = \frac{(t\ell)^{\gamma\ell}}{(\gamma\ell)!} 
\le \left(\frac{e\ell t}{\gamma\ell} \right)^{\gamma\ell} \le e^{-\gamma \ell}~,
\]
where we used Stirling's formula and the choice $t=\gamma/e^2$.
Putting everything together, we have that
\begin{eqnarray*}
\limsup_{n\to \infty} \PROB\left\{ \min_{\ell<i\le n} \psi(i)\le \max_{\ell \gamma/2 \le j \le \ell (1-\gamma/2)} \psi(j) \right\} 
 \le 
2\ell^2 \left(  e^{-\gamma \ell} + e^{-\gamma (\ell-1)/e^2}
 \right) \le \epsilon
\end{eqnarray*}
under our conditions for $\ell$, 
as desired.
\hfill $\square$

\subsection{Proof of Theorem \ref{thm:pathlower}}

Let $E$ be the event that either (1) vertex $i$ attaches to vertex $i-1$ for all $i=\ell+1,\ldots,2\ell$
or (2) vertex $\ell+1$ attaches to vertex $1$ and for all $i=\ell+2,\ldots,2\ell$, vertex
$i$ attaches to vertex $i-1$. On this event,
$T_{2\ell}$ is a path of $2\ell$ vertices such that the seed $P_{\ell}$
is on one of the two extremes of $T_{2\ell}$. The probability of this event is
$$
\dfrac{2}{\ell}\cdot \dfrac{1}{\ell+1} \cdot \cdots \cdot \dfrac{1}{2\ell-1} \geq 2\dfrac{\ell!}{(2\ell)!} \geq
2 (2\ell)^{-\ell}~.
$$
On this event, for $n\ge 2\ell$, for any seed-finding algorithm, the first and second halves of the path $T_{2\ell}$
are indistinguishable. At least one of the two halves of $T_{2\ell}$ is such that $H_n$ intersects that half in
at most $\ell/2$ vertices. Thus,
(conditionally on $E$), the algorithm misses at least half of the seed path, with probability $1/2$.
Hence
\[
 \PROB\left\{ |H_n \cap P_\ell| \le \frac{\ell}{2} \right\} \geq \frac{\PROB\{E\}}{2} \ge (2\ell)^{-\ell} \ge 
\epsilon
\]
whenever $\ell \leq  \frac{ \log(1/\epsilon)}{\log \log (1/\epsilon)}$ and $\epsilon\le e^{-e^2}$.

\subsection{Proof of Theorem \ref{thm:star}}

Let $k_\ell=(1+\gamma)\ell$. Again, we may assume that $k_\ell$ is an integer.
The seed finding algorithm we propose is slightly different. It is specifically tailored to the 
case when the seed tree to be found is a star. Let $v_n^* = \argmin_{i=1,\ldots,n} \psi(i)$ be the most
central vertex of $T_n$. We define $H_n$ as the set of vertices that includes $v_n^*$ and $k_\ell-1$ other
vertices $j$ with largest value of 
$\left|(T_n,v_n^*)_{j\downarrow}\right|$ among the neighbors of $v_n^*$ in $T_n$.
In other words, the algorithm outputs the most central vertex $v_n^*$ and those neighbors whose subtree
away from $v_n^*$ is largest.

First we recall that by Jog and Loh \cite[Theorem 4]{JoLo17a}, there exists a numerical constant $C$ such that,
if $\ell \ge C\log(1/\epsilon)$ and the uniform attachment tree is initialized with a star $E_\ell$ as seed of $\ell$
vertices and central vertex $1$, then 
\[
    \PROB\left\{v_n^*=1  \text{ for all $n=\ell+1,\ell+2,\ldots$ } \right\} \ge 1- \frac{\epsilon}{2}~, 
\]
that is, with probability at least $1- \epsilon/2$, the center of the seed star remains the most central vertex of
$T_n$ for all $n$.

Let $v_1\le v_2\le \cdots $ be the vertices that are attached to vertex $1$ (i.e., to the
center of the seed star $E_\ell$) in the uniform attachment process. (Thus, $v_1>\ell$.)
In view of the above-mentioned result of Jog and Loh, it suffices to show that for all $n$
sufficiently large, all vertices $v_j$ with  $j>\gamma\ell$  have $\left|(T_n,1)_{v_j\downarrow}\right|$
smaller than $\left|(T_n,1)_{i\downarrow}\right|$ for all vertices $i$
in the seed star $E_\ell$, with probability at least $1-\epsilon/2$. Thus, writing $g(i)=\left|(T_n,1)_{i\downarrow}\right|$,
we need to prove that
\begin{equation}
\label{eq:starcentr}
\limsup_{n\to \infty}   \PROB\left\{ \max_{j>\gamma\ell} g(v_j) < \min_{i=2,\ldots,\ell} g(i) \right\} > 1-\frac{\epsilon}{2}~.
\end{equation}

To prove (\ref{eq:starcentr}), first we write
\begin{equation}
\label{eq:star1}
   \PROB\left\{ \max_{j>\gamma\ell} g(v_j) \ge \min_{i=2,\ldots,\ell} g(i) \right\} 
\le \PROB\left\{ v_{\gamma\ell+1} \le m \right\} 
   + \PROB\left\{ \max_{v_j > m} g(v_j) \ge \min_{i=2,\ldots,\ell} g(i) \right\}~,
\end{equation}
where we take $m=\lfloor e^{\gamma\ell/4}\rfloor$. The first term on the right-hand side is the probability that more than $\gamma \ell$ vertices
are attached to vertex $1$ up to time $m$. In order to bound this probability, 
denote by $X_t$, for $t\ge \ell$, the number of vertices attached to vertex $1$ between time $\ell+1$ and $t$.
Thus, $X_\ell = 0$ and 
\[
\PROB\left\{ v_{\gamma\ell+1} \le m \right\}  = \PROB\left\{ X_m > \gamma\ell\right\}~.
\]
Since
\[
\EXP[X_t|X_{t-1}] = X_{t-1} + \frac{1}{t}~,
\]
\[
Y_t = X_t - \sum_{k=\ell+1}^{t} \frac{1}{k}, \ \ t \geq \ell +1
\]
is a martingale with respect to the filtration generated by $X_\ell, X_{\ell+1},\ldots$. 
Denote the corresponding martingale difference sequence by $Z_t = Y_t-Y_{t-1} = X_t-X_{t-1} - 1/t$.
By Markov's inequality,
\begin{equation}
\label{eq:markov}
\PROB \left\{ X_m > \gamma \ell \right\} 
 = 
\PROB\left\{ \sum_{j=\ell+1}^{m} Z_j + \sum_{j=\ell+1}^{m} \frac{1}{j} > \gamma \ell \right\} 
\leq  \frac{e^{ \sum_{j=\ell+1}^{m} \frac{1}{j}} \cdot\EXP\left[e^{ \sum_{j=\ell+1}^{m} Z_j } \right]}{e^{\gamma\ell}}~.
\end{equation}
In order to bound the right-hand side, observe that
\begin{eqnarray*}
 \EXP\left[e^{ Z_m}  | X_\ell,\ldots,X_{m-1} \right] 
&= & \EXP\left[e^{ X_m-X_{m-1} - \frac{1}{m}}  | \  X_\ell,\ldots,X_{m-1}\right] \\
&= & e^{ -X_{m-1} - \frac{1}{m}}  \EXP\left[e^{ X_m}  | \  X_\ell,\ldots,X_{m-1}  \right] \\
&= & e^{ -X_{m-1} - \frac{1}{m} }  \left( \frac{1}{m}e^{ X_{m-1}+1}+ \frac{(m-1)}{m}e^{X_{m-1}}\right) \\
&= & \frac{e^{ - \frac{1}{m} } }{m} ( e+ m-1)\\
& \leq & \frac{(m+2) e^{ - \frac{1}{m}}}{m}~,
\end{eqnarray*}
and therefore
\begin{eqnarray*}
\EXP\left[e^{ \sum_{j=\ell+1}^m Z_j } \right] 
&= &\EXP\left[ \EXP\left[e^{ \sum_{j=\ell+1}^m Z_j }  |  X_\ell,\ldots,X_{m-1} \right] \right] \\
&= & \EXP\left[e^{ \sum_{j=\ell+1}^{m-1} Z_j }  \EXP\left[e^{ Z_m}  | X_\ell,\ldots,X_{m-1} \right]\right] \\
&\leq &  \frac{(m+2) e^{ - \frac{1}{m} } }{m} \EXP\left[e^{ \sum_{j=\ell+1}^{m-1} Z_j } \right]~.
\end{eqnarray*}
Thus, by induction we obtain
\[
\EXP\left[e^{ \sum_{j=\ell+1}^m Z_j } \right] \leq \frac{(m+2)^2}{\ell^2} e^{ -\sum_{j=\ell+1}^m \frac{1}{j}}.
\]
Substituting into (\ref{eq:markov}), we get
\[
\PROB\left\{ v_{\gamma\ell+1} \le m \right\}  = \PROB\left\{ X_m > \gamma \ell\right\}  \leq \frac{(m+2)^2}{\ell^2 e^{\gamma\ell}}
  \le \frac{\epsilon}{4}
\]
by our choice of $m$ and by the condition on the value of $\ell$.
Hence, by (\ref{eq:star1}), it suffices to show that
\[
 \PROB\left\{ \max_{v_j> m} g(v_j) \ge \min_{i=2,\ldots,\ell} g(i) \right\} \le \frac{\epsilon}{4}~.
\]
We proceed by writing
\[
 \PROB\left\{ \max_{v_j> m} g(v_j) \ge \min_{i=2,\ldots,\ell} g(i) \right\} \le 
\sum_{i=2}^\ell \PROB\left\{ \max_{v_j> m} g(v_j)\ge g(i) \right\}~.
\]
Now fix $i\in \{2,\ldots,\ell\}$ and notice that $\max_{v_j> m} g(j)$ is bounded by the number of vertices $A$
attached to the tree formed by vertex $1$ and all vertices in the subtrees $(T_n,1)_{j\downarrow}$ for 
$j>m$ such that vertex $j$ is attached to vertex $1$.

Denoting $B=g(i)$
and $C=n-A-B$, note that, conditioned on the tree $T_m$, the triple $(A,B,C)$ behaves as the number of
red, blue, and white balls in a P\'olya urn in which initially (i.e., at time $m$) there is
one red ball, $B_m=\left|(T_m,1)_{i\downarrow}\right|$ blue balls, and
$m-1-\left|(T_m,1)_{i\downarrow}\right|$ white balls.
Hence, for each $i=2,\ldots,\ell$, we have
\begin{eqnarray*}
  \PROB\left\{ \max_{v_j> m} g(v_j)\ge g(i) \right\}
  & \le &
  \PROB\left\{A>B\right\}   \\
    & \le &
  \PROB\left\{A>B | B_m \ge \frac{m\epsilon}{32\ell^2} \right\} +   \PROB\left\{B_m < \frac{m\epsilon}{32\ell^2} \right\}~.
\end{eqnarray*}
In order to bound the second term on the right-hand side, note that
by the standard theory of P\'olya urns, $B_m$ has a beta-binomial distribution
with parameters $(m,1,\ell-1)$. Thus, $B_m$ is distributed as
a binomial random variable $\text{Bin}(m,\pi)$ where the parameter $\pi$ is
an independent $\text{Beta}(1,\ell-1)$ random variable.
Thus,
\begin{eqnarray*}
\lefteqn{
  \PROB\left\{B_m < \frac{m\epsilon}{32\ell^2} \right\}
} \\
  & \le &
  \PROB\left\{\text{Bin}(m,\epsilon/16\ell^2) < \frac{m\epsilon}{32\ell^2} \right\}
  + \PROB\left\{ \pi < \frac{\epsilon}{16\ell^2} \right\}
  \\
  & \le &
  e^{-m\epsilon/(128\ell^2)} + 1 -\left(1-\frac{m\epsilon}{16\ell^2}\right)^{\ell-1}
  \\
  & & \text{(by a standard binomial estimate and expressing the beta distribution)} \\
  & \le &
    e^{-m\epsilon/(128\ell^2)} + \frac{\epsilon}{16\ell} \\
    & & \text{(by the Bernoulli inequality)} \\
    & \le & \frac{\epsilon}{8\ell}
\end{eqnarray*}
whenever $\ell > (4\gamma)\left(\log(1/\epsilon)+\log\log(8\ell/\epsilon)
+ \log(128\ell^2)\right)$.
To finish the proof it remains to show that
\[
  \limsup_{n\to\infty}
  \PROB\left\{A>B | B_m \ge \frac{m\epsilon}{32\ell^2} \right\} \le
  \frac{\epsilon}{8\ell}~.
\]
But this follows from the fact that this limiting probability is
bounded by the the probability that a $\text{Beta}(1,m\epsilon/32\ell^2)$
random variable is greater than $1/2$ which is at most $2^{-m\epsilon/32\ell^2}$.
Since $m=\lfloor e^{\gamma\ell/4}\rfloor$, this is bounded by $\epsilon/(8\ell)$
for $\ell > (8/\gamma\vee C)\log(1/\epsilon)$, as desired.
 $\square$

\subsection{Proof of Theorem \ref{thm:oldtestament}}

Fix $\epsilon\in (0,1)$ and define $a=2\log(4/\epsilon)+1$ and $k_{\ell} = \frac{\ell}{3a}$.
A seed-finding algorithm with the desired property simply selects
the $k_\ell$ most central vertices. (Again, for simplicity of the presentation, we assume that
$k_\ell$ is an integer.)
With the notation introduced at
the beginning of this section, we define $H_n = H_{\psi}(k_{\ell})$.
We need to show that the $k_{\ell}$ most central vertices of $T_n$
are in $T_{\ell}$ with probability at least $1-\epsilon$
for all sufficiently large $n$.

The strategy of our proof is as follows. First we show that, with
probability at least $1-\epsilon/2$, the seed $T_\ell$ contains at least
$k_\ell$ ``deep'' vertices. Then we prove that
for all $n$ sufficiently large, all deep vertices of $T_\ell$ are
more central in $T_n$ than any vertex outside of the seed $T_\ell$.

We call a vertex $v\in T_\ell$ \emph{deep} if it has at least $a$
descendants, that is, if
\[
\left|(T_\ell,1)_{v\downarrow}\right| \ge a+1~.
\]
Denote by 
$\A_{\ell}$ the set of all deep vertices of $T_\ell$. 
Noticing that
\[
    \PROB\left\{ H_n \not\subset T_{\ell} \right\} \le
\PROB\left\{ |\A_{\ell}| \leq k_{\ell} \right\}  +
 \PROB\left\{ \exists v \in V(T_n)\backslash V(T_{\ell}), \exists u \in \A_{\ell} : \psi_n(v) \leq \psi_n(u) \right\}~,
\]
it suffices to show that 
\begin{equation}
\label{eq:oldtest1}
\PROB\left\{ |\A_{\ell}| \leq k_{\ell}
\right\} 
\le \frac{\epsilon}{2}~.
\end{equation}
 and
\begin{equation}
\label{eq:oldtest2}
\limsup_{n\to\infty}   \PROB\left\{ \exists v \in V(T_n)\backslash
  V(T_{\ell}), \exists u \in \A_{\ell} : \psi_n(v) \leq \psi_n(u)
\right\}
\le \frac{\epsilon}{2}~.
\end{equation}
(\ref{eq:oldtest1}) follows from inequality (\ref{eq:deepvertices}) in the
Appendix under the condition $\ell \ge 64a^2\log(22a/\epsilon)$.

It remains to prove (\ref{eq:oldtest2}). To this end, for $i\in \{1,\ldots,\ell\}$, denote
by $C_i$ the component of vertex $i$ in the forest obtained
by removing the edges of $T_{\ell}$ from $T_n$.
Then
\begin{eqnarray*}
  \PROB\left\{ \exists v \in V(T_n)\backslash
  V(T_{\ell}), \exists u \in \A_{\ell} : \psi(v) \leq \psi(u)
  |T_{\ell} \right\}
   \\
 \leq  \sum_{u \in \A_{\ell}} \sum_{k=1}^{\ell} \PROB\left\{ \exists v \in C_k\backslash\{k\}: \psi(v) \leq \psi(u) |T_\ell \right\}~.
\end{eqnarray*}
Now fix $T_\ell$ and vertices $k\in \{1,\ldots,\ell\}$ and $u \in \A_{\ell}$.
For any vertex $v \in C_k\backslash\{k\}$ such that $\psi(v) \leq \psi(u)$,
there are two possibilities:

\noindent
(1) either the largest subtree of $T_n$ rooted at $v$ is inside $C_k$,
in which case $|C_k| \geq \sum_{i \neq k}|C_i|$;

\noindent
(2) or the largest subtree of $T_n$ rooted at $v$ is  $\left(\bigcup_{i=1, i\neq k}^{\ell} C_i \right) \cup C_k'$ for some $C_k' \subset C_k$.
In this case, $\psi(v)\le \psi(u)$ implies that
\[
\sum_{i \in T_n \backslash (T_\ell,v)_{u\downarrow}} |C_i| \leq  |C_k|~.
\]
Since $u\in A_\ell$, this means that the left-hand side is dominated
by the number of red balls 
in a standard P\'olya urn with after $n-\ell$ draws initialized
with at least $a$ red, one blue, and $n-a-\ell-1$ white balls; while $|C_k|$
behaves like the number of blue balls in the same urn.

By the same calculations as in the proof of Theorem \ref{thm:path},
the probability of case (1) may be bounded by
\[
\limsup_{n\to\infty} \PROB\left\{ |C_k| \geq \sum_{i \neq k}|C_i| |T_\ell \right\}
= \limsup_{n\to\infty} \PROB\left\{ |C_k| \geq (n-\ell)/2 |T_\ell\right\}
 \le e^{-(\ell-1)/2} \le \frac{\epsilon}{4\ell^2}~.
\]
Similarly, the probability of case (2) satisfies
\[
\limsup_{n\to\infty} \PROB\left\{
\sum_{i \in T_n \backslash (T_\ell,v)_{u\downarrow}} |C_i| \leq  |C_k| |T_\ell\right\}
\le  e^{-(a-1)/2} \le \frac{\epsilon}{4\ell^2}
\]
by our choice $a=2\log(\ell^2/\epsilon)+1$.
This concludes the proof of (\ref{eq:oldtest2}) and hence
that of Theorem \ref{thm:oldtestament}.

\subsection{Proof of Theorem \ref{thm:oldtestamentlower}}

We prove the lower bound for the expected number of camouflaging vertices
by induction. To this end, fix a singleton $d$ and its parent $v$ in $T_\ell$.
For $j\ge \ell$, let
\[
E_j^{(v)} =\left\{\exists d'\in V(T_j)\backslash \{d\}: d'\sim v
\ \text{and} \ d',d \ \text{are leaves in} \ T_j\right\}~.
\]
Observe that $E_{2\ell}^{(v)}$ is the event that
$v$ is a camouflaging vertex. Consider the sequences
\begin{align*}
a_j &= \PROB\left\{E_j^{(v)}|T_\ell\right\} \\
c_j &= \PROB\left\{ \text{$d$ is a singleton in $T_j$}|T_\ell\right\}~.
\end{align*}
Now, observe that the event $E_{j+1}^{(v)}$ occurs if $E_j^{(v)}$ occurs
and the vertex $j+1$ is neither attached to $d$ nor to $d'$,
or if $d$ is a singleton of $T_{j}$ and the $j+1$ is attached to $v$. Thus
\[
a_{j+1} = a_j\cdot \left(1-\frac{2}{j}\right) + c_j\cdot \frac{1}{j}~.
\]
Multiplying both sides by $j(j-1)$, we get
\[
j(j-1)a_{j+1} = (j-1)(j-2)a_j + (j-1) c_j~.
\]
Summing over $j=\ell+1,\ldots,2\ell-1$,
\[
(2\ell-1)(2\ell-2) a_{2\ell} = \ell(\ell-1)a_{\ell+1} + \sum_{j=\ell+1}^{2\ell-1}(j-1)c_j~,
\]
which implies that
\[
a_{2\ell} \geq \frac{1}{(2\ell-1)(2\ell-2)} \sum_{j=\ell+1}^{2\ell-1}(j-1)c_j  \geq \frac{1}{4(\ell-1)} \sum_{j=\ell+1}^{2\ell-1} c_j~.
\]
Note that, for $j\in \{\ell+1,\ldots,2\ell-1\}$,
\begin{eqnarray*}
c_j & = & \prod_{k=\ell}^{j-1} \left(1-\frac{2}{k} \right)  \\
&\geq & \exp\left( -4\sum_{k=\ell}^{j-1} \frac{1}{k} \right)
\quad \text{(since $1-x \geq e^{-2x}$ for $x<3/4$)} \\
& \geq &
\exp\left(4 \log \ell - 4 \log j \right) \\
&> & \frac{\ell^4}{(2\ell)^4} =\frac{1}{16}~,
\end{eqnarray*}
and therefore
\begin{align*}
a_{2\ell} \ge \frac{1}{4(\ell-1)} \sum_{j=\ell+1}^{2\ell-1} c_j \geq \frac{1}{64}~.
\end{align*}
Let $P_\ell$ be the set of vertices in $T_\ell$ that are parents of a singleton. Then
\begin{align*}
\EXP[G_{\ell}|T_\ell] &= \EXP\left[ \sum_{v \in P_\ell} 1_{E_{2\ell}^{(v)}}|T_\ell \right] \\
&= \sum_{v \in P_\ell} \PROB\left\{E_{2\ell}^{(v)} |T_\ell \right\} \\
&\geq \frac{1}{64}|P_\ell|~,
\end{align*}
which implies that $\EXP G_{\ell}  \geq \frac{1}{64}\EXP |P_\ell|$.

It remains to bound the expected number of singletons $\EXP|P_\ell|$
in the uniform random recursive tree $T_\ell$.
Write $S_k=|P_k|$ and note that $S_k$ 
equals the number of parents of singletons in $T_k$.

When a new vertex is attached to the tree $T_k$, we lose one
singleton if the new vertex is attached to the parent of a singleton.
This happens with probability $S_k/k$.
If a the new vertex is attached to a singleton,
then the number remains the same.
If the new vertex is attached to some vertex that is not a
leaf nor a parent of a singleton, then, the number of singletons
also remains unchanged.
Finally, if the new vertex is attached to a leaf
that is not a singleton, the number of singletons increases by $1$.
Thus, denoting the number of leaves of $T_k$ by $L_k$,
\begin{align*}
  \EXP[S_{k+1}|T_k] &=
  (S_k-1)\frac{S_k}{k} + S_k \left(\frac{S_k}{k}+1-\frac{S_k}{k}
  -\frac{L_k}{k}\right) + (S_k+1)\left(\frac{L_k}{k}-\frac{S_k}{k}\right) \\
& = \left(1-\frac{2}{k} \right)S_k + \frac{L_k}{k}~.
\end{align*}
Taking expectations and using the fact that
$\EXP L_k=k/2$, we have that $\EXP S_\ell=\ell/6$.
Summarizing, the expected number of camouflaging vertices satisfies
\[
\EXP G_{\ell} \geq \frac{1}{64}\cdot \frac{\ell}{6}
=\frac{\ell}{384}~.
\]
We prove the second inequality of Theorem \ref{thm:oldtestamentlower}
using the \emph{bounded differences inequality} of McDiarmid \cite{McD89}
(see also \cite[Theorem 6.2]{BoLuMa13}).

Observe that given $T_\ell$, there is a bijection between the set of
recursive trees of size $2\ell$ containing $T_\ell$ as subgraph and the
set $\mathcal{S} = [\ell] \times \cdots \times [2\ell-1]$. The bijection
is simply given by 
associating the vector $\kappa = (a_{\ell+1},\cdots, a_{2\ell})$ to the
recursive tree $T(\kappa)$ where the vertex $k \in [\ell+1,2\ell]$ is attached
to the vertex $a_{k}$, starting by $T_\ell$ until obtaining
$T_{2\ell}$. 
Then we may consider the set
$\mathcal{S}$ as the set of recursive trees with $2\ell$ vertices that contain
$T_\ell$ as subtree.

Importantly, the components of $\kappa$ that represent 
the uniform random recursive tree $T_{2\ell}$ are independent
random variables.

Given $T_\ell$, consider the function $g: \mathcal{S} \to \R$ such that
$g(T_{2\ell})$ is the number of camouflaging vertices.

By the bounded differences inequality, it suffices to
show that, given $T,T' \in \mathcal{S}$,
if $T$ and $T'$ differ by exactly one coordinate, then
$|g(T) - g(T')| \leq 2.$

To this end, let $v \in V(T_n)$ be a parent of a singleton $d$. 
$v$ is a camouflaging vertex of a tree $T=(a_{\ell+1},\cdots,a_{2\ell})$ if and only if
\begin{enumerate}
\item $d \notin \{a_{\ell+1},\cdots,a_{2\ell}\}$;

\item $\exists k \in \{\ell+1,\cdots, 2\ell\}\backslash \{a_{k+1},\cdots,a_{2\ell}\}$ such that $a_k=v$.
\end{enumerate}

Now, consider $T=(a_{\ell+1},\cdots,a_{2\ell})$,
$T'=(b_{\ell+1},\cdots,b_{2\ell})$ two trees with $a_r \neq b_r$ for
some $r$ and $a_j=b_j$ for $j \neq r$. For a camouflaging vertex $v$
in $T$ (with corresponding singleton $d$ in $T_\ell$) not to be a
camouflaging vertex in $T'$, 
it is necessary (but not sufficient) that either
\begin{enumerate}
\item $b_r$ is a child of $v$, 
\item or $a_r=v$~.
\end{enumerate}
Similarly, for a not camouflaging vertex $v$ in $T$ (with
corresponding singleton $d$ in $T_\ell$), to be a camouflaging vertex in $T'$ it is necessary that either
\begin{enumerate}
\item $a_r$ is a descendant of $v$, 
\item or $b_r=v$~.
\end{enumerate}
Thus, $|g(T)-g(T')|\leq 2$, and the bounded differences condition
is satisfied, proving the second inequality of Theorem \ref{thm:oldtestamentlower}.

\section{Appendix}

Devroye \cite{Dev91} proved a central limit theorem for the number of
vertices with $k$ descendants in a uniform random recursive tree.
In particular, if $L_{k,n}$ denotes the
the number of
vertices with $k$ descendants in a uniform random recursive tree of $n>k+1$
vertices, then Devroye shows that
\[
\EXP L_{k,n} = \frac{n-k-1}{(k+1)(k+2)} + \frac{1}{k+1}
      = \frac{n+1}{(k+1)(k+2)} 
\]
and, for any fixed $k$, as $n\to\infty$,
\[
   \frac{L_{k,n}- \frac{n}{(k+1)(k+2)}}{\sqrt{n\sigma_k^2}}
   \]
   converges, in distribution, to a standard normal random variable,
   where
\[
\sigma_k^2 = \frac{1}{(k+1)(k+2)}\left(1-\frac{1}{(k+1)(k+2)}\right)
   -\frac{2}{(k+1)(k+2)^2}+ \frac{1}{(k+1)^2(2k+3)}~.
\]
Devroye's proof is based on representing $L_{k,n}$ as a sum of $(k+1)$-dependent
indicator random variables and on a central limit theorem of Hoeffding
and Robbins \cite{HoRo48} for such sums. In this paper we need a non-asymptotic
version of Devroye's theorem. Quantitative, Berry-Esseen-type versions
of the Hoeffding-Robbins limit theorem are available via Stein's method,
see, for example, Rinott \cite[Theorem 2.2]{Rin94}. On the other hand, a simple
bound may be proved by combining Devroye's representation with a concentration
inequality of Janson \cite[Corollary 2.4]{Jan04} for sums of
dependent random variables,
to obtain the following:

\begin{proposition}
If $L_{k,n}$ denotes the
the number of
vertices with $k$ descendants in a uniform random recursive tree of $n>k+1$,
then for all $t>0$,
\[
\PROB\left\{ L_{k,n} \ge \EXP L_{k,n} + t \right\}
\le \exp\left(\frac{-8t^2(k+2)}{25(n+(k+1)(k+2)t/3}\right)
 \]
and   
\[
\PROB\left\{ L_{k,n} \le \EXP L_{k,n} - t \right\}
\le \exp\left(\frac{-8t^2(k+2)}{25n}\right)~.
 \]
\end{proposition}

Note that the number of vertices with \emph{at least} $k$ descendants
$M_{k,n}=\sum_{i=k}^{n-1} L_{i,n}=n-\sum_{i=0}^{k-1} L_{i,n}$ has expected value
\[
\EXP M_{k,n}= \EXP  \sum_{i=k}^{n-1} L_{i,n}=n-\sum_{i=0}^{k-1} \EXP L_{i,n}
= \frac{n+1}{k+1} -1~,
\]
and therefore
\begin{eqnarray*}
  \PROB\left\{ M_{k,n} \le \frac{n+1}{k+1} -1 -t \right\} 
& = & \PROB\left\{ \sum_{i=0}^{k-1} L_{i,n} \ge \sum_{i=0}^{k-1} \EXP
  L_{i,n}  +t\right\}  \\
& \le & \sum_{i=0}^{k-1} \PROB\left\{ L_{i,n} \ge \EXP  L_{i,n}
  +\frac{t}{k} \right\} \\
& \le & k \exp\left(\frac{-8t^2}{25k(n+(k+1)t/3}\right)~.
\end{eqnarray*}
In particular, by generously bounding constants, we get
\begin{equation}
\label{eq:deepvertices}
 \PROB\left\{ M_{k,n} \le \frac{n}{3k} \right\} \le
 k\exp\left(- \frac{1}{32} \frac{n}{k^2}\right)~.
\end{equation}

\medskip
\noindent
{\bf Acknowledgements.} We thank Luc Devroye, Mikl\'os R\'acz, and
Tommy Reddad for interesting conversations on the topic of the paper.

\bibliographystyle{plain}
\bibliography{findingtheseed}

\end{document}